\documentclass{article}

\usepackage[affil-it]{authblk}
\usepackage[top=1in, bottom=1in, left=0.75in, right=0.75in]{geometry}
\usepackage{graphicx}
\usepackage{amssymb,amsfonts,amsmath}
\usepackage{placeins}
\usepackage{setspace}

\newcommand{\eq}[1]{Eq.~\eqref{eq:#1}}
\newcommand{\spc}{\mbox{ }}
\newcommand{\url}[1]{#1}

\newcommand{\FigIa}{1(a)}
\newcommand{\FigId}{1(d)}

\usepackage{color}
\usepackage{ulem}

\begin{document}
\doublespacing
\normalem

\title{The influence of societal individualism on a century of tobacco use: modelling the prevalence of smoking}

\author[1]{John C. Lang}
\author[2]{Daniel M. Abrams}
\author[1]{Hans De Sterck}
\affil[1]{Department of Applied Mathematics, University of Waterloo}
\affil[2]{Department of Engineering Sciences and Applied Mathematics \& Northwestern Institute on Complex Systems \& Department of Physics and Astronomy, Northwestern University}

\maketitle


\begin{abstract}
	Smoking of tobacco is predicted to cause approximately six million deaths worldwide in 2014. Responding effectively to this epidemic requires a thorough understanding of how smoking behaviour is transmitted and modified.  Here, we present a new mathematical model of the social dynamics that cause cigarette smoking to spread in a population. Our model predicts that more individualistic societies will show faster adoption and cessation of smoking.  Evidence from a new century-long composite data set on smoking prevalence in 25 countries supports the model, with direct implications for public health interventions around the world. Our results suggest that differences in culture between societies can measurably affect the temporal dynamics of a social spreading process, and that these effects can be understood via a quantitative mathematical model matched to observations.
\end{abstract}

\section{Introduction and Motivation}
	In the fifty years since the first report of the Surgeon General's Advisory Committee on Smoking and Health \cite{SurgGen64} the smoking epidemic has been responsible for more than 20 million deaths in the United States alone \cite{SurgGen14, Warner81}, and continues to be responsible for over 6 million deaths worldwide each year \cite{Jha09, MathersLoncar06}. The strong social component of the dynamics of smoking prevalence has been modelled mathematically \cite{CavanaTobias08, SharomiGumel08, LevyBauerLee06, RoweEtAl92}, and examined statistically through analysis of social network data \cite{ChristakisFowler08} and survey data \cite{HoskingEtAl09, UngerEtAl01, Janis83}. However, whereas previous works tend to focus on the micro-level, in this paper we investigate how social aspects of smoking affect its prevalence at the societal level. 

	Significant inter-country differences exist in smoking prevalence \cite{Pampel10}. For example, Fig.~\FigIa\spc shows smoking prevalence estimates over most of the past century for Sweden and the USA, obtained from surveys and cigarette consumption data (see Section 3 and Appendix \ref{A:Data}). In both countries, smoking prevalence increased rapidly starting from the early decades of the 20\textsuperscript{th} century and reached a peak in the 1960s--1980s era when the adverse health effects of smoking became widely known \cite{SurgGen64}, after which smoking prevalence declined rapidly. However, there are conspicuous differences between the curves: the rate of smoking adoption and cessation before and after the peak is much greater in the US than in Sweden, and the peak in prevalence in the US occurs much earlier than in Sweden.

\begin{figure}
	\centering
        	\includegraphics{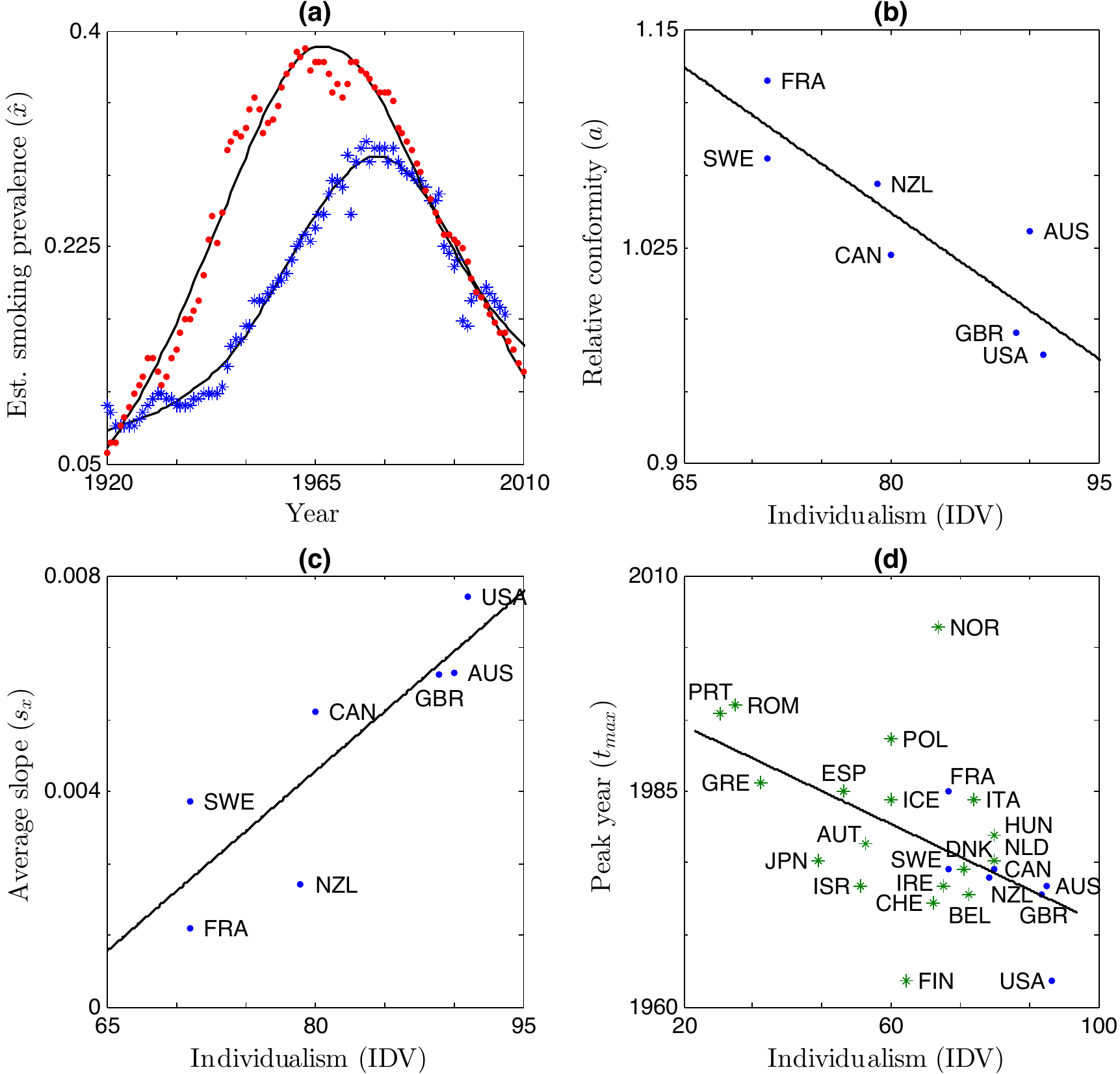}
        	\caption{Model results and relation with tobacco consumption and individualism data. (a) Output of \eq{1} (solid line) fitted to historical tobacco use data for the United States (dots) and Sweden (asterisks); US IDV is 91, Swedish IDV is 71. (b) Relative conformity parameter $a$ versus IDV ($\rho = -0.87$, $p = 0.011$) for seven OECD countries. (c) Average slope $s_x$ versus IDV ($\rho = 0.85$, $p = 0.015$) for seven OECD countries. (d) Peak year $t_{max}$ in cigarette consumption versus IDV in a set of 25 countries ($\rho = -0.524$, $p = 0.008$). The seven OECD countries considered for the mathematical model (and displayed in panels (b-c)) are indicated by dots ($\rho = -0.76$, $p = 0.047$), and the remaining 18 countries are indicated by asterisks. For panels (b-d) the line of best fit is given by a solid line. For panel (d) the line of best fit is calculated using data from all 25 countries.}
	\label{fig:1}
\end{figure}

\begin{figure}
	\centering
        	\includegraphics[width=8cm]{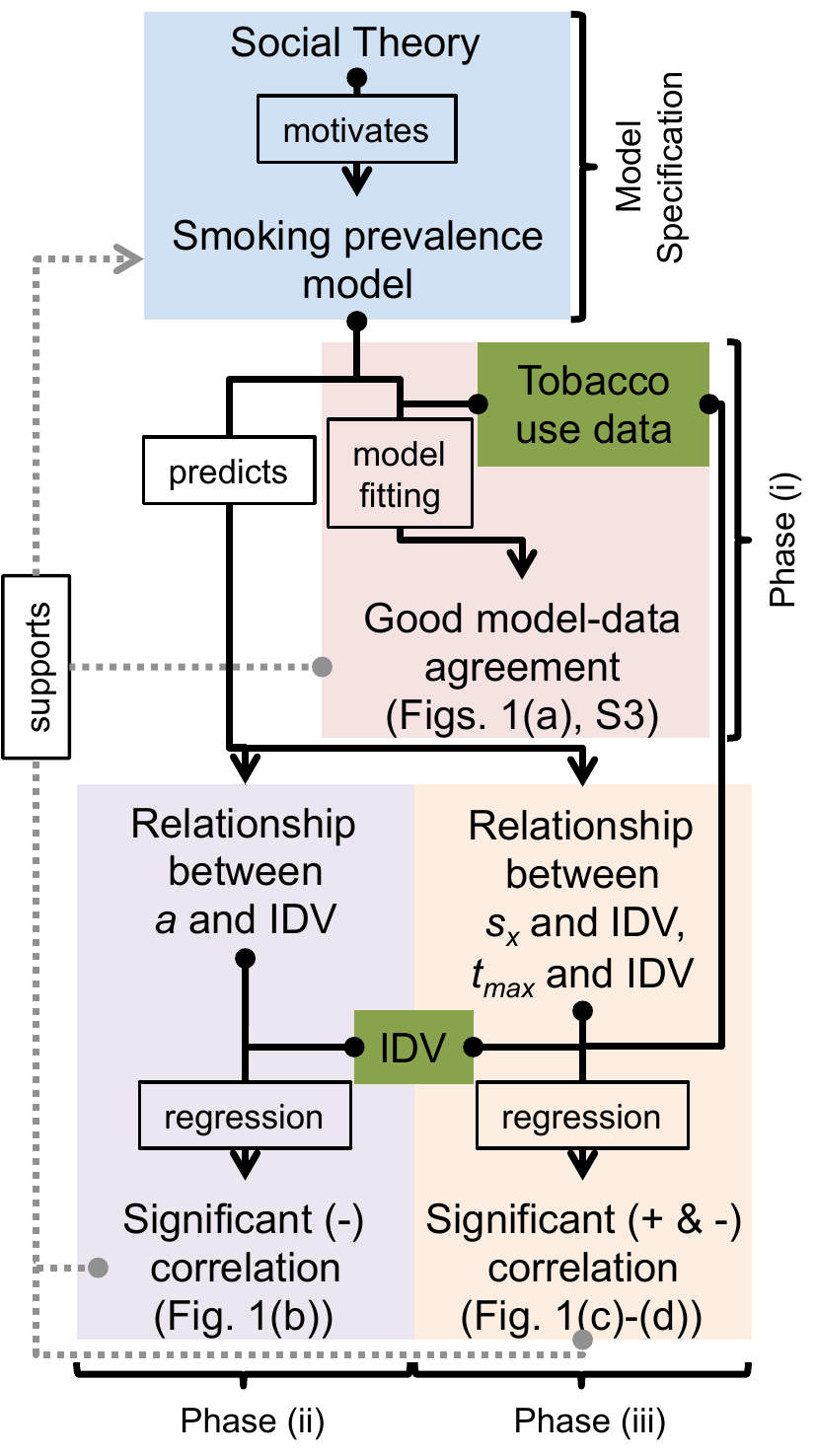}
        	\caption{Model specification and testing of model predictions/interpretation in three phases.}
        	\label{fig:flow}
\end{figure}

	In this paper we present a new mathematical model of the social spreading of smoking that incorporates the concepts of individual utility (including awareness of health effects), peer influence and social inertia. We propose an interpretation for our model in the context of societal individualism/collectivism and test the model's predictions and interpretation in three separate phases, see Fig.~\ref{fig:flow}. First, we compile smoking prevalence data spanning the past century for seven OECD countries and find good agreement between this data and the fitted model. Second, the parameter in our model that we interpret as the degree of societal collectivism (see Section \ref{Sec:Spec}), and that we fit to smoking prevalence data, is found to be significantly negatively correlated to an established measure of societal individualism (Hofstede's IDV \cite{HofstedeEtAl10}). Thirdly, the central role played by societal individualism/collectivism in our model motivates us to investigate directly the role that individualism (as measured by Hofstede's IDV) plays in 
observed historical tobacco use data. Specifically, we find that IDV is significantly positively correlated to the average rate of increase in smoking prevalence ($s_x$) in seven OECD countries for which smoking prevalence data is available, and that it is significantly negatively correlated to the peak year of tobacco consumption ($t_{max}$) for 25 countries in which tobacco consumption data are available. These findings are interpreted according to our modelling framework, which offers an explanation for the compelling effects of individualism/collectivism on smoking prevalence.

\section{Model Specification}
\label{Sec:Spec}
	We begin formulating our model by observing that individuals derive utility from smoking via two mechanisms. First, they derive utility directly from the act of smoking (\emph{individual utility}). Second, they derive utility from social interaction with other smokers (\emph{social utility}). We note that social utility commonly manifests itself in the form of peer influence or peer pressure \cite{CalavoArmengol10, Centola05}. We then proceed using a modelling framework that explicitly accounts for the effect of competition between individual and social utilities, and that was first applied to explore the temporal dynamics of language death and religious affiliation as binary choice problems \cite{Abrams03, Abrams11}. Specifically, we propose the model 
\begin{equation}
	\label{eq:1}
	\frac{dx}{dt} = b \spc [(1-x) \spc x^a u_x - x \spc (1-x)^a (1-u_x)],
\end{equation}
where $x\in [0,1]$ is the fraction of smokers in the population (i.e., the prevalence), $u_x \in [0,1]$ is the individual utility from smoking, and the constant $b > 0$ determines the timescale of the equation. The interpretation of the positive term in \eq{1} is, therefore, that non-smokers $1-x$ take up smoking at a rate proportional to the \emph{total utility} derived from smoking, $x^au_x$, which is the weighted product of the individual utility from smoking $u_x$ and the social utility from interactions with other smokers $x$, with weighting determined by the constant parameter $a$. Since societies with large $a$ weigh social utility more heavily than individual utility when calculating total utility, we call $a$ the \emph{relative conformity parameter}. We therefore interpret societies with large $a$ to be more collectivistic (or less individualistic) than societies with small $a$. The interpretation of the negative term in \eq{1}, which models smoking cessation, follows analogously: smokers $x$ cease smoking at a rate proportional to the total utility derived from non-smoking, $(1-x)^a(1-u_x)$, which is the weighted product of the individual utility from non-smoking $u_y = 1-u_x$ and the social utility from interactions with other non-smokers $1-x$, where we have normalized individual utilities from smoking $u_x$ and from non-smoking $u_y$ such that $u_x+u_y=1$.  We note that this modelling framework is conceptually consistent with the findings presented in \cite{HoskingEtAl09}: that personal attitudes about smoking have a stronger influence on smoking behaviour in individualistic countries than in collectivistic countries. 

	Next, we observe that a combination of factors, including advances in our understanding of the health effects of smoking and public policy initiatives designed to curb smoking, have likely reduced individual utility from smoking ($u_x$) over the past century. Thus, in a significant departure from previous work that treats individual utility as a constant \cite{Abrams03, Abrams11}, we account for this decline in individual utility by using the cumulative number of scholarly articles on the health effects of smoking ($n(t)$) as a proxy for the reduction in individual utility over the past century. Specifically, following the principle of temporal discounting \cite{Frederick02}, we assume that each additional article published is discounted by the factor $\delta \in [0,1]$ so that for year $t$
\begin{equation}
	\label{eq:2}
	u_x(t) = u_\infty + \delta^{n(t)} \spc (u_0 - u_\infty),
\end{equation}
where $u_0$ and $u_\infty$ are the limiting individual utilities from smoking when there is no knowledge and perfect knowledge of the adverse effects of smoking, respectively. Here, $u_0$, $u_\infty$ and $\delta$ are parameters to be fitted to observational data. Fig.~\ref{fig:S1} shows the cumulative number of scholarly articles on the health effects of smoking $n(t)$ (see Appendix \ref{A:Proxy} and ESM data file \emph{lang\_supporting\_data1.txt}), and the resulting individual utility from smoking $u_x(t)$, as defined in \eq{2}, for various values of $\delta$. We remark that this approach leads to better fits between model output and observational data than alternatives that do not directly take into account the effect of increased scientific understanding of health effects (see Appendices \ref{A:Fit} and \ref{A:Alt}).

\section{Data}
	We note that \eq{1} subject to \eq{2} requires the fitting of four parameters per country ($x_0=x(t_0)$, $a$, $u_0$, and $u_{\infty}$) and two parameters $b$ and $\delta$ that we take equal for all countries in the data set (see Appendix \ref{A:Fit}). We determine these parameters by fitting them to historical smoking prevalence data. Such historical data (mostly obtained by surveys) is available for a set of 24 OECD countries, but unfortunately is limited to, on average, only 21.5 observations over a period of 31.4 years spanning 1960--2012 \cite{OECD}. As such, it misses much of the crucial period in the earlier parts of the 20th century during which smoking steadily gained popularity in many countries. However, historical national cigarette consumption data is available for the same 24 OECD countries plus Romania for an average of 78.4 observations over a period of 82.2 years spanning 1900-2012 \cite{ForeyEtAl, ForeyEtAl02}. In order to exploit the more abundant consumption data in our prevalence model, we assume a linear relationship between national cigarette consumption and smoking prevalence for all times where both measurements are available. We then apply this relationship to the full cigarette consumption data set to produce an estimate of historical smoking prevalence ($\hat{x}$) over nearly the entire past century. Sufficient data for this procedure is available in seven OECD countries (for additional details, see Appendix \ref{A:Data}): Australia, Canada, France, New Zealand, Sweden, the United Kingdom, and the United States.

\section{Results: Testing the Model}
  \subsection{Phase (i): Direct test}
   Figure \FigIa\spc shows the fit of our model to data sets from the United States and Sweden (additional fits and parameter values are displayed in Fig.~\ref{fig:S3} and Table \ref{tab:S4}).  The good agreement that we found with all data sets provides support for the model.  

  \subsection{Phase(ii): Test of model implications for $a$}
   If the model and its interpretation are correct, then we expect that the fitted ``relative conformity parameter'' $a$ will capture something meaningful about the individualism/collectivism of a society.  To test this, we compare with Hofstede's IDV, an established metric for societal individualism \cite{HofstedeEtAl10} that has been evaluated in most countries. Panel (b) of Fig.~\ref{fig:1} shows the comparison.  As expected, the relative conformity parameter $a$ shows significant negative correlation with Hofstede's IDV (negative because $a$ increases with collectivism while IDV decreases with it).  This concordance with independently assessed individualism values supports our model.
   
  \subsection{Phase (iii): Test of model implications for slope and peak year}
   Besides the correlation of $a$ with collectivism, we note that another prediction is implicit in model \eqref{eq:1}. As the relative conformity parameter increases, the model requires that changes in smoking prevalence occur more slowly (this is true for the range of $a$ and $u$ values corresponding to the observational data).  Put another way, societies with higher levels of individualism should experience faster changes in smoking prevalence. Intuitively, when smoking prevalence is low the lack of existing smokers inhibits smoking initiation more strongly in a collectivistic society than in an individualistic society. Thus, we expect the average rate of increase in a collectivistic society to be smaller than in an individualistic society. In contrast, when smoking prevalence is high, and once the deleterious health effects of smoking become widely known and negatively impact individual utility from smoking, the presence of existing smokers inhibits smoking cessation more strongly in a collectivistic society than in an individualistic society.  In both cases collectivism acts as a break on change in the status quo (higher cultural inertia \cite{Zarate12, Henrich03}). Panel (c) of Fig.~\ref{fig:1} demonstrates that this is indeed the case: the average slope $s_x$ of the smoking prevalence curves leading up to the peak (see \eq{sx} of Appendix \ref{A:sx}) increases with Hofstede's IDV (shown) and decreases with $a$ (see Fig.~\ref{fig:S4}(a)).  Correlations are significant (see Table \ref{tab:S1}).
	
	This reasoning further suggests that the peak year for smoking prevalence $t_{max}$ should be later in collectivistic societies and earlier in individualistic countries. As shown in Fig.~\FigId\spc, the raw observational data are consistent with this prediction: $t_{max}$ is significantly negatively correlated with IDV (shown) and significantly positively correlated with $a$ (see Fig.~\ref{fig:S4}(b)). Note that our assumption of a linear relationship between national cigarette consumption and smoking prevalence is not needed to establish $t_{max}$, so Fig.~\FigId\spc is independent of any model assumptions. 

\section{Discussion and Conclusion}
    We have proposed a quantitative mathematical model of the social spreading of smoking that is derived from basic principles well-documented in the sociology and social psychology literature. The model appears to match real-world smoking prevalence data from a variety of countries well (to our knowledge, the largest historical data set of this type ever compiled), and all predictions of the model appear to be supported by the data. In particular, the model predicts that the level of individualism or collectivism of a society may significantly affect the temporal dynamics of smoking prevalence: the strong influence of the personal utility of smoking (and its decrease due to increased awareness of adverse health effects) is predicted to lead to faster adoption and cessation of smoking in individualistic societies than in more collectivistic societies.

    It has previously been argued that social support mechanisms in collectivistic societies make it more likely that a person will stop smoking \cite{Sun04, Triandis88} based on findings that social support (supportive counselors) can help people to adhere to decisions to quit smoking \cite{Janis83}. We find that, to the contrary, cessation of smoking occurs more slowly in collectivistic societies.   Our model suggests that this is so because social inertia will inhibit decisions to stop smoking more strongly in collectivistic societies than in individualistic societies. 
    
    These results have significant implications for combating the ongoing smoking epidemic. For example, they imply that interventions designed to discourage smoking should be tailored differently in societies or social groups whose cultures differ in how they value individualism versus collectivism \cite{Pentecostes99}. More broadly, these results demonstrate that differences in culture can measurably affect the dynamics of a social spreading process, and that a mathematical model can help to illuminate this phenomenon.

Despite the good match between model predictions and data, a number of limitations remain.  We have made an implicit ``mean-field'' approximation in taking social utility to be a function of the overall smoking prevalence $x$, rather than the local smoking prevalence among contacts in an individual's social network.  Similarly, we have taken individual utility to be uniform across the population (though not in time), whereas a more detailed model might allow for individual variation.  As a mild justification for these assumptions, we point out that analysis of a similar model in another context \cite{Abrams11} suggests that inclusion of more detail will not change qualitative predictions.

    We claim that the correlation of individualism with faster societal change results from a causative influence as predicted by our model. Other factors such as income levels also correlate with individualism, and it's possible that what we observe is ultimately also related to GDP or other variables.  We certainly cannot exclude that there may be other causative factors. For example, our model in its current form is incapable of explaining differences in smoking prevalence between genders and why these inter-gender differences vary between countries \cite{Pampel10}. Nevertheless, we remark that many previously proposed causative factors for differences in observed inter-country smoking dynamics can be accounted for within our modelling framework. In particular, beliefs about the harmful effects of smoking, the price of cigarettes, socioeconomic status and inequality, and government regulation have all been cited as potential factors affecting the differences observed in inter-country smoking dynamics \cite{Pampel10}. Each of these factors can be interpreted within our modelling framework. For example, beliefs about the harmful effects of smoking, as well as the price of cigarettes, both likely contribute directly to individual utility derived from smoking ($u_x$) and from non-smoking ($u_y$). Moreover, socioeconomic status may affect individual utility from smoking indirectly by affecting an individual's tolerance for risk and/or how they discount future rewards and costs (i.e. how they discount their future health status) \cite{Griskevicius11}. Addressing the model's inability to account for gender differences in smoking prevalence and quantifying the relationship between other causative factors and model parameters are potential areas for future work.

    We also welcome future work comparing a variety of social contagion phenomena across societies.  Our model suggests that the increased cultural inertia in collectivistic societies would lead to slower change across a wide spectrum of spreading processes (those where important changes occur in personal utility), a hypothesis that could be supported or rejected by further study.

\section*{Acknowledgments}
	We would like to thank Dr.~James Fowler for insightful remarks on an early draft of this manuscript. DMA thanks the James S.~McDonnell Foundation for support through grant \#220020230. JL and HDS acknowledge support from NSERC of Canada. 

\nocite{OECDPop}

\appendix
\setcounter{figure}{0}
\setcounter{equation}{0}

\renewcommand\thefigure{\thesection.\arabic{figure}}
\renewcommand\theequation{\thesection.\arabic{equation}}
\renewcommand\thetable{\thesection.\arabic{table}}

\section{Appendix - Materials and Methods}
\subsection{Smoking prevalence and cigarette consumption data}
\label{A:Data}
	We consider smoking prevalence $x(t) \in [0,1]$ for 24 OECD countries which we download from the OECD iLibrary online statistical database \cite{OECD} in Excel format. We also consider manufactured cigarette consumption (in grams) per person per day $c(t)$ for the same 24 OECD countries plus Romania (which is a non-OECD country) \cite{ForeyEtAl,ForeyEtAl02}. When available, cigarette consumption data is downloaded directly from the International Smoking Statistics (Web Edition) website \cite{ForeyEtAl} in Excel format. Cigarette consumption data for countries not included in the  International Smoking Statistics (Web Edition) are retrieved from the International Smoking Statistics (2\textsuperscript{nd} Ed.) \cite{ForeyEtAl02} by manually transferring these entries into Excel. We make these data available in CSV format in the ESM data file \emph{lang\_supplementary\_data1.txt}, which contains four columns: country number as it appears in Table~\ref{tab:S2}, year ($t$), measurement ($x(t)$ or $c(t)$), and type of measurement (0 indicates a smoking prevalence measurement, while 1 indicates a cigarette consumption measurement).
\begin{table}[hp]
	\centering
	\caption{Summary of data on smoking prevalence $x$ and cigarette consumption (in grams) per person per day $c$}
	\begin{tabular}{lllccccc}
		\hline
		&&& \multicolumn{2}{c}{} & \multicolumn{3}{c}{Cigarette consumption per person}\\
		&&& \multicolumn{2}{c}{Smoking prevalence ($x$)} & \multicolumn{3}{c}{per day ($c$)}\\
		No. 	& Country & Abbrev. & Obs. Period& No. of Obs. & Obs. Period & No. of Obs.  & Source\\ \hline
		1	& Australia 		& AUS 	& 1964--2010 	& 16 			& 1920--2010 	& 91 	& \cite{ForeyEtAl}\\
		2	& Austria		& AUT 	& 1972--2006	& 5			& 1923--2004	& 82	& \cite{ForeyEtAl}\\
		3	& Belgium		& BEL		& 1997--2008	& 4			& 1921--2011	& 91	& \cite{ForeyEtAl}\\
		4	& Canada 		& CAN 	& 1964--2011 	& 29 			& 1920--2010 	& 91	& \cite{ForeyEtAl}\\
		5	& Denmark		& DNK		& 1970--2010	& 41			& 1920--2010	& 91	& \cite{ForeyEtAl}\\
		6	& Finland		& FIN		& 1978--2011	& 34			& 1920--2009	& 90	& \cite{ForeyEtAl}\\
		7	& France 		& FRA 	& 1960--2010 	& 22 			& 1900--2010 	& 93	& \cite{ForeyEtAl}\\ 
		8	& Greece		& GRE		& 1998--2009	& 6 			& 1920--1995	& 76 	& \cite{ForeyEtAl02}\\
		9	& Hungary		& HUN		& 1994--2009 	& 4			& 1920--2012	& 87	& \cite{ForeyEtAl}\\
		10	& Iceland 		& ICE		& 1987--2012 	& 26 			& 1932--1995 	& 64	& \cite{ForeyEtAl02}\\
		11	& Ireland 		& IRE		& 1973--2007 	& 14 			& 1920--1995 	& 76	& \cite{ForeyEtAl02}\\
		12	& Israel		& ISR		& 1996--2010	& 8 			& 1967--1995	& 29	& \cite{ForeyEtAl02}\\
		13	& Italy		& ITA		& 1980--2012	& 23			& 1905--2010	& 73	& \cite{ForeyEtAl}\\
		14	& Japan		& JPN		& 1965--2011	& 47			& 1920--2007	& 88	& \cite{ForeyEtAl}\\
		15	& Netherlands	& NLD		& 1966--2011	& 39			& 1923--1995	& 67	& \cite{ForeyEtAl02}\\
		16	& New Zealand 	& NZL		& 1976--2012 	& 28 			& 1920--2009 	& 90	& \cite{ForeyEtAl}\\
		17	& Norway		& NOR	& 1973--2012	& 40			& 1927--2011	& 85	& \cite{ForeyEtAl}\\
		18	& Poland		& POL		& 1996--2009	& 4			& 1925--1995	& 43	& \cite{ForeyEtAl02}\\
		19	& Portugal		& PRT		& 1987--2006	& 4			& 1940--1995	& 56	& \cite{ForeyEtAl02}\\
		20	& Romania		& ROM	& -- 			& 0			& 1920--1995	& 52	& \cite{ForeyEtAl02}\\
		21	& Spain		& SPA 	& 1985--2011	& 11			& 1920--2010	& 87	& \cite{ForeyEtAl}\\
		22	& Sweden 		& SWE 	&1980--2011 	& 32 			& 1920--2006 	& 87	& \cite{ForeyEtAl}\\
		23	& Switzerland	& CHE		& 1992--2007	& 4			& 1934--2009	& 76	& \cite{ForeyEtAl}\\
		24	& United Kingdom & GBR		& 1960--2010 	& 38 			& 1905--2009 	& 105& \cite{ForeyEtAl}\\
		25	& United States 	& USA		& 1965--2011 	& 36 			& 1920--2010 	& 91	& \cite{ForeyEtAl}\\\hline
	\end{tabular}
	\label{tab:S2}
\end{table}

	Recall that the cigarette consumption data are more dense and span a larger time interval than smoking prevalence data, see Table~\ref{tab:S2}. Thus, since our model is specified in terms of smoking prevalence, we estimate smoking prevalence from cigarette consumption in order to exploit the much richer cigarette consumption data for model fitting purposes. First, we assume a linear relationship between smoking prevalence $x(t)$ and smoking consumption $c(t)$ 
$$
	x(t) = C c(t) + B.
$$

\begin{figure}[hp]
  \centering
  \includegraphics[width=0.9\linewidth]{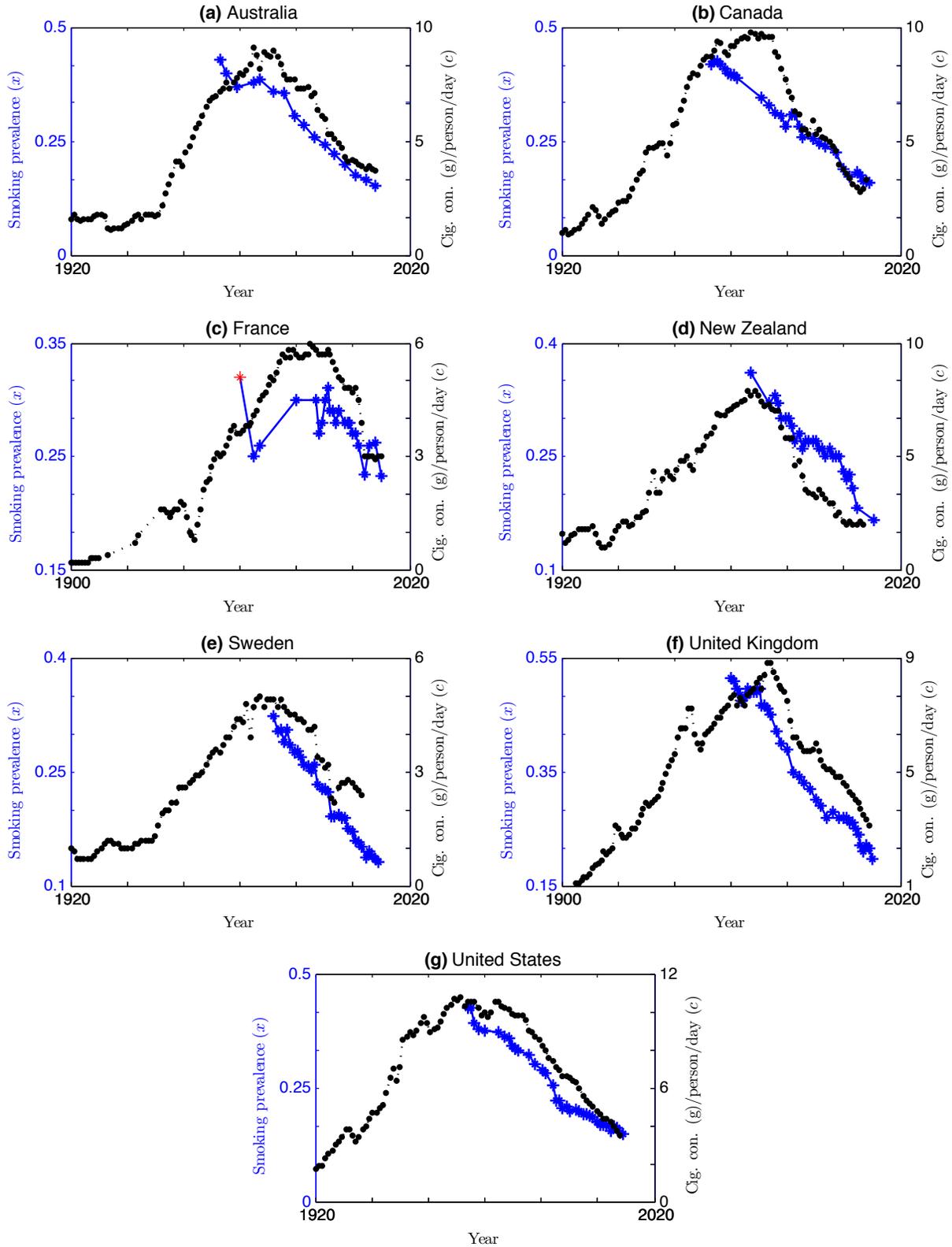}
  \caption{Raw smoking prevalence ($x$) and raw cigarette consumption ($c$) data. (Left axis - blue asterisks) Raw smoking prevalence data and (right axis - black dots) raw cigarette consumption data (in grams per person per day). A single outlier for smoking prevalence ($x$) for the country of France (panel c) is denoted with a red asterisk.}
	\label{fig:S2}
\end{figure}

	Next, we calculate estimates $\widehat{C}$ and $\widehat{B}$ by regressing smoking prevalence $x(t)$ on tobacco consumption $c(t)$ for all years for which both measurements are available. The results of this regression are summarized in Table~\ref{tab:S3}, which illustrates that the assumption that $x$ and $c$ are linearly related does not hold equally well for all countries. In order to restrict ourselves to the cases where the assumption of linearity between $x$ and $c$ is valid we restrict ourselves to the seven OECD countries with $R^2 \geq 0.7$, $p<0.001$, and $n_{obs}\geq 15$: Australia, Canada, France, New Zealand, Sweden, the United Kingdom, and the United States. We display the raw data for these seven OECD nations in Fig.~\ref{fig:S2}. The smoking prevalence for these seven OECD countries is then estimated from tobacco consumption using the relationship
$$
	\hat{x}(t) = \widehat{C} c(t) + \widehat{B}.
$$
\begin{table}[hp]
	\centering
	\caption{Estimates $\widehat{C}$, $\widehat{B}$}
	\begin{tabular}{llllll}
		\hline
		Country & $\widehat{C}\times 10^2$ & $\widehat{B}\times 10^2$ & $R^2$ & $p$ & $n_{obs}$\\ \hline
		Australia 		& $4.5 \pm 1.3$	& $-0.3 \pm 8.8$ 	& 0.80 	& $3.2\times10^{-6}$ 	& 16\\
		Austria		& $0.0 \pm 4.9$ 	& $24.2 \pm 32.4$	& 0.00	& 0.99			& 4\\
		Belgium		& $2.6\pm20.3$	& $13.0 \pm 81.5$	& 0.13	& 0.64			& 4\\
		Canada 		& $3.5 \pm 0.5$ 	& $6.3 \pm 3.8$ 		& 0.87	& $3.0\times10^{-13}$	& 28\\
		Denmark		& $0.0 \pm 9.2$	& $40.5 \pm 44.4$	& 0.00	& 0.99			& 41\\
		Finland		& $2.0 \pm 0.7$	& $15.8 \pm 2.8$	& 0.55	& $1.0\times10^{-6}$	& 32\\
		France		& $1.8 \pm 0.5$ 	& $19.1 \pm 2.5$ 	& 0.72	& $6.3\times10^{-7}$ 	& 22\\ 
		Greece		& --			& --				& -- 		& --				& 0\\
		Hungary		& $1.9 \pm 1.6$	& $17.4 \pm 11.2$ 	& 0.93	& $3.5\times10^{-2}$	& 4\\
		Iceland 		& $4.9 \pm 1.2$	& $0.9 \pm 7.0$ 		& 0.93 	& $2.6\times10^{-5}$ 	& 9\\
		Ireland 		& $5.4 \pm 1.1$	& $-4.0 \pm 7.4$ 	& 0.93 	& $1.7\times10^{-6}$ 	& 11\\
		Israel		& --			& --				& -- 		& -- 				& --\\
		Italy			& $4.8 \pm 2.5$	& $-0.3 \pm 13.2$	& 0.47	& $6.1\times 10^{-4}$	& 21\\
		Japan		& $1.3 \pm 3.2$	& $25.7 \pm 27.2$	& 0.02	& $0.43$			& 43\\
		Netherlands	& $4.8 \pm 3.2$	& $20.5 \pm 15.0$	& 0.32	& $4.7\times10^{-3}$	& 23\\
		New Zealand 	& $2.0 \pm 0.3$	& $18.8 \pm 1.4$ 	& 0.86 	& $2.6\times10^{-12}$	& 27\\
		Norway		& $-7.2 \pm 4.3$& $50.1 \pm 10.6$	& 0.24	& $1.6\times10^{-3}$	& 39\\
		Poland		& --			& --				& --		& -- 				& 0\\
		Portugal		& --			& --				& --		& --				& 1\\
		Romania		& --			& -- 				& --		& -- 				& 0\\
		Spain		& $6.0 \pm 6.2$ 	& $-7.4 \pm 41.7$	& 0.38	& $5.7\times10^{-2}$	& 10\\
		Sweden 		& $5.4 \pm 0.6$ 	& $4.3 \pm 2.3$ 		& 0.92 	& $1.7\times10^{-15}$	& 27\\
		Switzerland	& $2.8 \pm 5.6$	& $7.2 \pm 38.6$	& 0.69	& 0.17			& 4\\
		United Kingdom & $5.6 \pm 0.7$	& $1.6 \pm 4.5$ 		& 0.88	& $5.3\times10^{-18}$	& 37\\
		United States 	& $3.6 \pm 0.3$	& $-0.1 \pm 2.3$ 	& 0.95 	& $1.1\times10^{-22}$	& 35\\\hline
	\end{tabular}\\
	\raggedright
	$\pm$ indicates 95\% confidence intervals. We report $R^2$ values for the linear regression of $x$ on $c$, the $p$-value of the correlation between $x$ and $c$, and the number of years for which both $x$ and $c$ measurements are available, $n_obs$.
	\label{tab:S3}
\end{table}

	Note that survey-based prevalence data are susceptible to noise stemming from variations in the survey methodology. In particular, prior to performing the linear regression of $x$ on $c$ for France, we removed the outlier $x(1960)=0.32$ since it is inconsistent with the rest of the data for France, see Fig.~\ref{fig:S2}(c). Specifically, the Grubbs test on $x / \hat{x}$ indicates that the 1960 data point is a significant outlier ($p < 0.05$). This can also be seen intuitively: from $t=1960$ until the next measurement at $t=1965$ smoking prevalence drops from $x(1960)=0.32$ to $x(1965)=0.25$ (a decrease of 21.9\%), while cigarette consumption steadily increases from $c(1960)=3.6$ to $c(1965)=4.1$ (an increase of 13.9\%). Given the population in France in 1960 (45.5 million) and in 1965 (48.6 million) \cite{OECDPop}, this would correspond to an increase in the average mass of cigarettes smoked (in grams) per smoker per day from 11.3 to 16.4 (an increase of 45.1\%) over a short 5 year period. This is in sharp contrast with the relatively stable relationship between $x$ and $c$ for France's remaining data points and justifies the exclusion of the outlier $x(1960)=0.32$. With the outlier removed, France satisfies our data quality requirements for inclusion in the set of seven OECD countries ($R^2\geq0.7$, $p<0.001$, and $n_{obs}\geq15$).

\subsection{Proxy data $n(t)$: articles published on the health effects of smoking}
\label{A:Proxy}
\begin{figure}[h]
  \centering
  \includegraphics[width=0.9\linewidth]{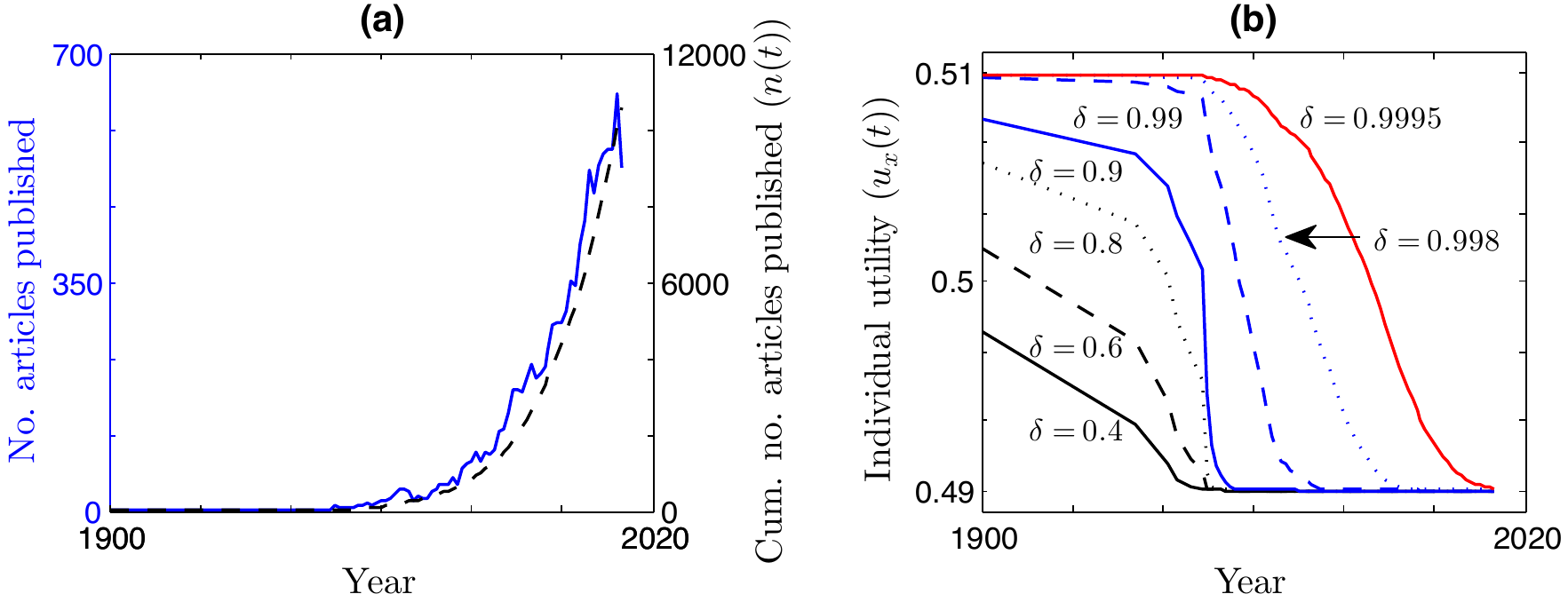}
  \caption{Articles published on the health effects of smoking. Articles retrieved by Scopus with search terms (i)-(iv) and individual utility profiles from \eq{2} for varying values of $\delta$. (a) (Left axis - blue, solid) Number of articles published per year and (right axis - black, dashed) cumulative number of articles published $n(t)$. (b) Discounted utility $u_x(t)$ from \eq{2} with $u_0 = 0.51$ and $u_\infty = 0.49$, using cumulative number of articles published $n(t)$. (Solid black) $\delta = 0.4$, (dashed black) $\delta = 0.6$, (dotted black) $\delta = 0.8$, (solid blue) $\delta = 0.9$, (dashed blue) $\delta = 0.99$, (dotted blue) $\delta = 0.998$, and (solid red) $\delta = 0.9995$.}
	\label{fig:S1}
\end{figure}

	In order to implement \eq{1} subject to \eq{2} we require data on the cumulative number of articles published on the health effects of smoking $n(t)$. We calculate $n(t)$ by performing a search of the online research database Scopus for papers with
\begin{enumerate}
	\renewcommand{\labelenumi}{(\roman{enumi})}
	\item tobacco, smok*, or cigar* in the title, and
	\item death, illness, mortality, risk*, tumour*, tumor*, or cancer in the title, and
	\item medicine, dentistry, nursing, veterinary, health professions, or multidisciplinary in the subject area, and
	\item plant*, mosaic, botany, smog, fog, and soot not in the title.
\end{enumerate}	 
Items (i)-(iii) are search terms included in order to select for papers researching the health effects of smoking, whereas items (iv) are search terms excluded in order to prevent selection of papers researching the tobacco mosaic virus (plant*, mosaic, botany) and the health effects of atmospheric smoke (smog, fog, soot). This provides us with $n(t)$ for integer $t$, where time $t$ is measured in years. We make the article data available in CSV format in the ESM data file, \emph{lang\_supplementary\_data2.txt}, which contains three columns: year ($t$), number of articles published in year $t$, and cumulative number of articles published up to and including year $t$ ($n(t)$). To calculate $n(t)$ for non-integer and missing values of $t$ we use linear interpolation, see Fig.~\ref{fig:S1}(a). Furthermore, Fig.~\ref{fig:S1}(b) displays $u_x(t)$ from \eq{2} using $n(t)$ calculated above for various discount factors $\delta$ and with $u_0=0.51$ and $u_{\infty}=0.49$. (For comparison, see Table~\ref{tab:S4}\spc for model-fitted values of $\delta$, $u_0$ and $u_{\infty}$.)

\subsection{Model fitting}
\label{A:Fit}
	We fit \eq{1} to the estimated prevalence, $\hat{x}(t)$. To reduce the dimensionality of the optimization problem, we assume that certain \emph{universal parameters} are constant across countries. Specifically, we assume that $b$ and $\delta$ are universal parameters, and that $x_i(t_{i,0})=x_{i,0}$, $a_i$, $u_{i,0}$, and $u_{i,\infty}$ are \emph{local parameters} for country $i$, where $t_{i,0}$ is the first year for which cigarette consumption data ($c$), and hence estimated smoking prevalence data ($\hat{x}$), are available. We denote the smoking prevalence estimated above for country $i$ at time $t$ by $\hat{x}_i(t)$. The time series of estimated smoking prevalences for country $i$ is then denoted by the vector $\widehat{X}_i$. Analogously, we denote the time series of smoking prevalences predicted by \eq{1} for country $i$ by $\widetilde{X}_i$. We solve \eq{1} using the Matlab differential equation solver \emph{ode45}.

	Using the Matlab function \emph{lsqcurvefit} we now proceed as follows:
\begin{enumerate}
	\item Holding universal parameters constant, for each country $i$ we find the $x_{i,0}$, $a_i$, $u_{i,0}$, and $u_{i,\infty}$ that minimize $E_i = \|\widetilde{X}_i-\widehat{X}_i\|_2^2$.
	\item Holding local parameters constant for each country $i$, we find the $b$ and $\delta$ that minimize $E = \sum_i\|\widetilde{X}_i-\widehat{X}_i\|_2^2$.
	\item Repeat steps (1) and (2) until either
	\begin{enumerate}
		\item the change in the objective function $\sum_i\|\widetilde{X}_i-\widehat{X}_i\|^2$ is below tolerance $tol$, or
		\item the number of iterations exceeds a limit $max_{itn}$.
	\end{enumerate}
\end{enumerate}

\begin{figure}[hp] 
  \centering
  \includegraphics[width=0.9\linewidth]{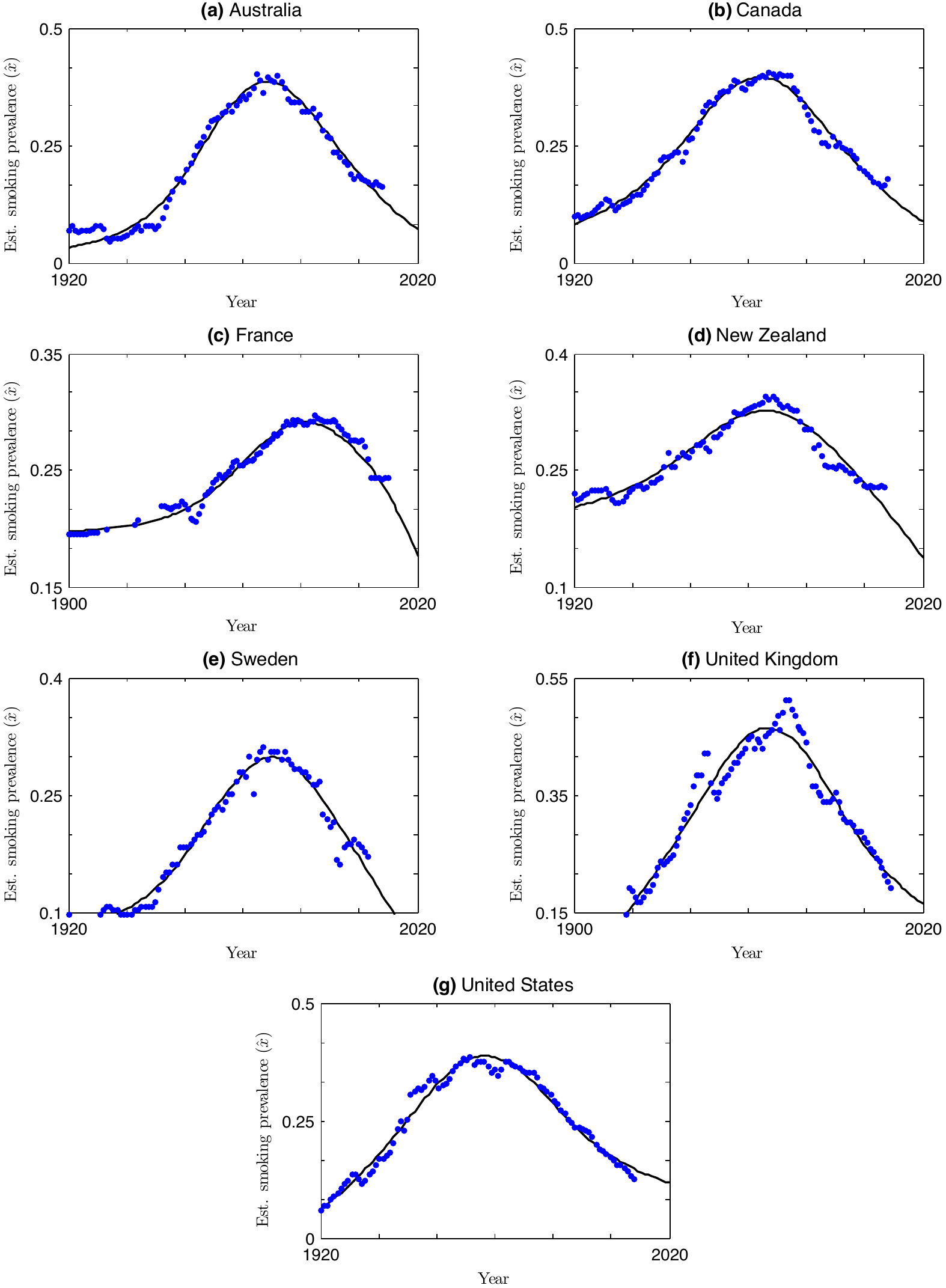}
  \caption{The result of fitting \eq{1} to the estimated smoking prevalence $\hat{x}$. Estimated smoking prevalence values $\hat{x}$ are given by blue dots.}
  \label{fig:S3}
\end{figure}

	We perform the optimization with the initial guess $u_{i,0}\equiv0.51$, $u_{i,\infty}\equiv0.49$, $x_{i,0} = \hat{x}_i(t_{i,0})$, $a_i=1$, $b=1$, and $\delta=0.9985$. We also provide the optimization algorithm \emph{lsqcurvefit} with constraints
\begin{align*}
	0&\leq a_i, b \leq 2 \mbox{ and}\\
	0&\leq x_{i,0}, u_{i,0}, u_{i,\infty}, \delta \leq 1,
\end{align*}
and with parameters $tol = 10^{-6}$ and $max_{itn} = 150$. The fitting procedure terminates after 114 iterations, the results of which are recorded in Table~\ref{tab:S4} and Fig.~\ref{fig:S3}. Complete model simulation code with all necessary data files are available upon request.
\begin{table}[h]
	\centering
	\caption{The result of fitting Model \eq{1} to the estimated smoking prevalence $\hat{x}$}
	\begin{tabular}{ccc} 
		\hline
	\multicolumn{3}{c}{Universal parameters and Total Error ($E$)}\\
		\hspace{8mm}$b$\hspace{8mm} & \hspace{8mm}$\delta$\hspace{8mm} & $E$ \\ \hline
		1.049 & 0.9981 & 0.163\\
		\hline
	\end{tabular}\\
	\vspace{10mm}
	\begin{tabular}{lccccc}
		\hline
	 	Country & \multicolumn{5}{c}{Local parameters and local error ($E_i$)}\\ 
			($i$) 		& $a_i$	& $x_{i,0}$ 	& $u_{i,0}$ & $u_{i,\infty}$	& $E_i$\\ \hline
		Australia		& 1.035	& 0.033 	& 0.551 	& 0.484 	& 0.032\\
		Canada 		& 1.020	& 0.083	& 0.530 	& 0.483	& 0.020\\
		France 		& 1.121 	& 0.198 	& 0.543 	& 0.524 	& 0.004\\
		New Zealand		& 1.062 	& 0.202 	& 0.525 	& 0.504 	& 0.012\\
		Sweden 		& 1.076 	& 0.077 	& 0.555 	& 0.503 	& 0.015\\
		United Kingdom	& 0.976 	& 0.079 	& 0.513 	& 0.478 	& 0.060\\
		United States	& 0.963 	& 0.063 	& 0.513 	& 0.470 	& 0.024\\\hline
	\end{tabular}
	\label{tab:S4}
\end{table}

\subsection{Alternative to temporal discounting of individual utility}
\label{A:Alt}
	We have found in more detailed analysis that using the discounting utility formula of \eq{2} produces a better fit (significantly lower total error $E$) than alternatives, e.g. constant utility $u_{i,x}(t)\equiv u_{i,x}$ or step-function utility 
$$
	u_{i,x}(t) = \left\{\begin{array}{ll} u_{i,0} & \mbox{if } t<t^*_i \\ u_{i,\infty} & \mbox{if }t\geq t^*_i \end{array}\right.,
$$
where $t^*_i$ is a threshold parameter whose value is determined by the fitting procedure. This is consistent with the expectation that increasing knowledge of health effects has indeed influenced the individual utility from smoking over the past century.

\subsection{Calculation of average slope}
\label{A:sx}
\begin{figure}[h] 
  \centering
  \includegraphics[width=0.9\linewidth]{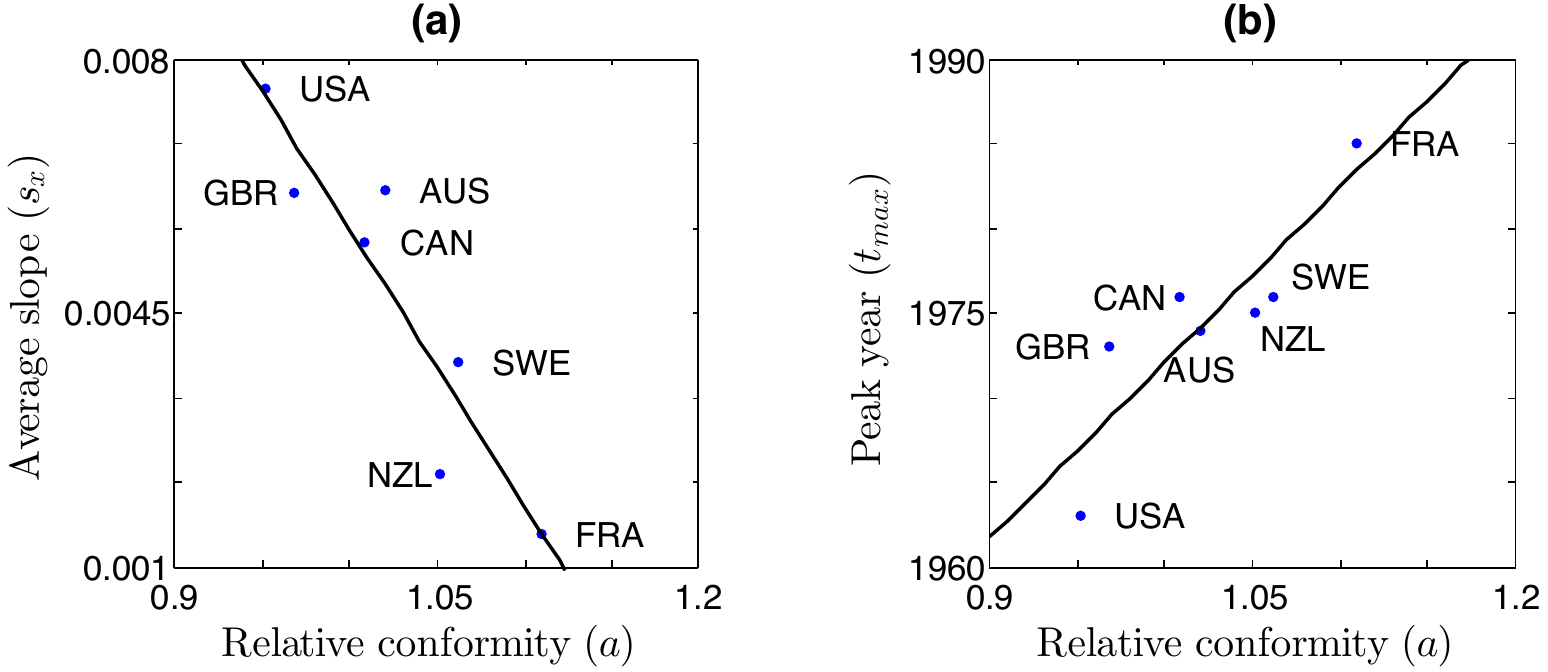}
  \caption{Average slope $s_x$ and peak year of smoking prevalence $t_{max}$ versus relative conformity parameter $a$ for seven OECD countries. (a) Average slope $s_x$ versus relative conformity parameter $a$ ($\rho = -0.92$, $p=0.003$). (b) Peak year $t_{max}$ versus relative conformity parameter $a$ ($\rho=0.88$,  $p=0.009$). The line of best fit is given by a solid line.}
  \label{fig:S4}
\end{figure}

We estimate the average slope of the estimated prevalence versus time curve in the period leading up to its peak as:
\begin{equation}
	\label{eq:sx}
	s_x = \frac{\hat{x}(t_{max}) - \hat{x}(t_0)}{t_{max}-t_0},
\end{equation}
where $t_0 = 1920$ is the first year for which smoking prevalence estimates are available in the subset of seven OECD countries, and where $t_{max}$ are recorded in Table \ref{tab:S5}. Correlation between fitted relative conformity parameter $a$ and average slope $s_x$ (and peak year $t_{max}$) is displayed in Fig.~\ref{fig:S4} and recorded in Table \ref{tab:S1}. 

\begin{table}[hp]
	\centering
	\caption{Hofstede's Individualism Index IDV and peak year in cigarette consumption ($t_{max}$)}
	\begin{tabular}{lcc}
		\hline
		Country \hspace{15mm}$\mbox{ }$ & IDV & Peak year ($t_{max}$) \\ \hline
		Australia 		& 90		& 1974\\
		Austria		& 55		& 1979\\
		Belgium		& 75		& 1973\\
		Canada 		& 80		& 1976\\
		Denmark		& 74		& 1976\\
		Finland		& 63		&1963\\
		France		& 71		& 1985\\ 
		Greece		& 35		& 1986\\
		Hungary		& 80		& 1980\\
		Iceland 		& 60 		& 1984\\
		Ireland 		& 70		& 1974\\
		Israel		& 54		& 1974\\
		Italy			& 76		& 1984\\
		Japan		& 46		& 1977\\
		Netherlands	& 80		& 1977\\
		New Zealand 	& 79		& 1975\\
		Norway		& 69		& 2004\\
		Poland		& 60		& 1991\\
		Portugal		& 27		& 1994\\
		Romania		& 30		& 1995\\
		Spain		& 51		& 1985\\
		Sweden 		& 71		& 1976\\
		Switzerland	& 68		& 1972\\
		United Kingdom & 89		& 1973\\
		United States 	& 91		& 1963\\ \hline
	\end{tabular}
	\label{tab:S5}
\end{table}
\begin{table}[h]
	\centering
	\caption{Individualism index IDV and relative conformity $a$ are significantly correlated with and $a$, $s_x$, and $t_{max}$}
	\label{tab:1}
 	\begin{tabular}{@{\vrule height 10.5pt depth4pt  width0pt}lcccc}
		\hline
 		& \multicolumn{3}{c}{7-country subset} & 25-country set\\
		& $a$		& $s_x$		& $t_{max}$	& $t_{max}$\\\hline
	IDV	& -0.87 (0.011)	& 0.85 (0.015)	& -0.76 (0.047)	& -0.53 (0.006)\\
	$a$	& -- 			& -0.92 (0.003)	& 0.88 (0.009)	& --\\\hline
 	\end{tabular}\\ 
	\raggedright
	Correlations between IDV, $a$, $s_x$, and $t_{max}$ are recorded for the seven-country subset. Correlation between IDV and $t_{max}$ is recorded for the full set of 25 countries. $p$-values are in parentheses. All correlations are significant at the 95\% confidence level.
	\label{tab:S1}
 \end{table}

\bibliographystyle{vancouver}

\end{document}